\documentclass{birkjour}

\usepackage{amssymb}
\usepackage{amsmath}
\usepackage{mathtools}
\usepackage{float}
\usepackage[colorlinks=true,pdftex,unicode=true,linktocpage,bookmarksopen,hypertexnames=false]{hyperref}
\usepackage{tikz-cd}
\usepackage{caption}

\usepackage{tikz}
\usetikzlibrary{braids}

\newcommand{\Ann}{\operatorname{Ann}}
\newcommand{\Soc}{\operatorname{Soc}}
\newcommand{\Aut}{\operatorname{Aut}}

\newcommand{\id}{\operatorname{id}}
\newcommand{\Ret}{\operatorname{Ret}}

\makeatletter
\numberwithin{equation}{section}
\numberwithin{figure}{section}
\numberwithin{table}{section}

\newtheorem{thm}{Theorem}[section]

\theoremstyle{plain}
\newtheorem{lem}[thm]{Lemma}

\newtheorem{pro}[thm]{Proposition}
\newtheorem{cor}[thm]{Corollary}
\newtheorem{exa}[thm]{Example}
\newtheorem{problem}[thm]{Problem}
\newtheorem{defn}[thm]{Definition}
\newtheorem{convention}[thm]{Convention}

\makeatother

\title{Skew braces: a brief survey}
\author{Leandro Vendramin}

\address{
Department of Mathematics, Vrije Universiteit Brussel, Pleinlaan 2, 1050 Brussel, Belgium}
\email{Leandro.Vendramin@vub.be}

\begin{document}

\begin{abstract}
  Our primary focus is on the theory of skew braces, specifically exploring their connection with combinatorial solutions to the Yang--Baxter equation. Skew braces have recently emerged as intriguing algebraic structures, and their link to the Yang--Baxter equation adds further depth and significance to their study.
  Throughout this article, we place particular emphasis on various problems and conjectures that arise within this field. These open questions serve to stimulate further research and investigation, as we strive to enhance our understanding of skew braces and their role in the study of the Yang--Baxter equation.
\end{abstract}

\maketitle

\tableofcontents

\section{Introduction}

The Yang--Baxter equation (YBE) holds significant importance in both pure mathematics and physics, serving as a fundamental concept in these fields. Its origins can be traced back to Yang's research on statistical mechanics, and later Baxter employed the YBE in his solution to the 8-vertex model. Since then, the YBE has been extensively studied, leading to exciting developments in various areas of mathematics and physics. Connections to knot theory, braid theory, operator theory, Hopf algebras, quantum groups, 3-manifolds, and the monodromy of differential equations have been established.

One of the fundamental challenges related to the YBE is the quest for solutions. Interestingly, a fascinating family of solutions can be constructed through linear extensions of maps $r\colon X\times X\to X\times X$, where $X$ is a basis of a vector space. These solutions, known as \emph{combinatorial solutions}, capture the combinatorial aspects of the YBE.

To construct and understand combinatorial solutions with specific properties, algebraic structures that act upon the solutions play a crucial role. These structures provide a framework for studying and characterizing the properties of solutions. Consequently, the algebraic structures underlying the YBE have become a major focus of research.

Among the various algebraic structures, skew braces have emerged as a particularly notable concept. Inspired by previous works \cite{MR3177933,MR2278047}, skew braces were introduced in \cite{MR3647970} and have found applications in diverse areas of mathematics. They have been connected to
triply factorizable groups \cite{MR4427114}, Garside theory~\cite{MR2764830,MR3374524}, regular subgroups~\cite{MR2273982,MR2486886}, ring theory~\cite{MR3824447,CJKVA}, flat manifolds~\cite{MR3291816}, pre-Lie algebras~\cite{Bai2023,MR4484785}, Rota--Baxter operators~\cite{MR4370524}, and Hopf--Galois structures~\cite{MR3763907,ST}.

In recent times, the theory of skew braces has also found applications in the realm of physics~\cite{MR4294149,MR4399582}. These applications demonstrate the interdisciplinary nature and broad relevance of skew braces beyond pure mathematics.

The main objective of this work is to provide a comprehensive overview of the fundamental theory of skew braces. We discuss their key properties, with a specific focus on their connection to solutions of the YBE. Additionally, we present a selection of  open problems that challenge researchers in the field. These problems encompass both the algebraic structure of skew braces and their associated solutions. By addressing these open problems, we aim to contribute to the advancement of research in this exciting field.

\section{Definitions and examples}

\begin{defn}
  \label{def:brace}
  A \emph{skew brace} is a triple $(A,+,\circ)$, where $(A,+)$ and $(A,\circ)$
  are (not necessarily abelian)
  groups and
  \begin{equation}
    \label{eq:compatibility}
    a\circ(b+c)=a\circ b-a+a\circ c
  \end{equation}
  holds for all $a,b,c\in A$. The groups
  $(A,+)$ and $(A,\circ)$ are respectively
  the \emph{additive} and \emph{multiplicative} group
  of the skew brace $A$.
\end{defn}

In~\eqref{eq:compatibility}, it is customary to grant multiplication higher precedence than addition. As a consequence of~\eqref{eq:compatibility}, we observe that the neutral elements of $(A,+)$ and $(A,\circ)$ coincide. This shared neutral element will be denoted as $0$ or $0_A$. In the context of a skew brace, the inverse of an element $a$ with respect to the circle operation is denoted by $a'$.
From \eqref{eq:compatibility}, the
following formulas
\[
  a\circ(-b+c)=a-a\circ b+a\circ c, \quad
  a\circ(b-c)=a\circ b-a\circ c+a
\]
hold for all $a,b,c$.

\begin{defn}
  Let $\mathcal{X}$ be a family of groups. A skew brace is said to be
  of \emph{$\mathcal{X}$-type} if its additive group belongs to $\mathcal{X}$.
\end{defn}

A significant and notable family of skew braces is that of \emph{skew braces of abelian type}. These skew braces possess an abelian additive group, and were initially introduced by Rump in~\cite{MR2278047}. They have been a subject of considerable interest and study. Another important family of skew braces is that of \emph{skew braces of nilpotent type}, which presents its own distinct properties and characteristics \cite{MR3957824}.

\begin{exa}
  \label{exa:trivial}
  If we consider an additive group $A$, it can be regarded as a skew brace by defining the operation $a\circ b = a + b$ for all $a, b \in A$. These skew braces are called \emph{trivial skew braces}.
\end{exa}

\begin{exa}
  If $A$ is an additive group, we can define $a\circ b = b + a$ for all $a, b \in A$. Then $A$ becomes a skew brace. These skew braces are called \emph{almost trivial skew braces}. 
\end{exa}

\begin{exa}
  \label{exa:WX}
  Let $A$ be an additive group
  and $B$ and $C$ be subgroups of $A$ such that $B\cap C=\{ 0\}$ and $A=B+C$.
  In this case, one says that $A$ admits an {\em exact factorization} through the subgroups $B$ and $C$.  Each $a\in A$ can be written uniquely as $a=b+c$, for some $b\in B$ and $c\in C$.  The map
  \[
    f\colon B\times C^{\mathrm{op}}\to A,\quad
    (b,c)\mapsto b+c,
  \]
  is bijective. Thus
  $A$ is a group with the operation
  \begin{equation*}
    %\label{eq:factorization}
    a\circ a_1=f\left(f^{-1}(a)f^{-1}(a_1)\right),\qquad a,a_1\in A.
  \end{equation*}
  Write $a=b+c$ and
  $a_1=b_1+c_1$ for some $b,b_1\in B$ and $c,c_1\in C$. A direct calculation shows that
  \begin{align*}
    a\circ a_1=b+a_1+c.
  \end{align*}
  Hence $(A,\circ)$ is a group
  isomorphic to $B\times C^{\mathrm{op}}\cong B\times C$. Moreover, routine calculations
  show that $(A,+,\circ)$ is a skew brace.
\end{exa}

\begin{lem}
  \label{lem:lambda}
  Let $A$ be a skew brace and $a\in A$. The map
  \[
    \lambda\colon (A,\circ)\to\Aut(A,+), \quad
    a\mapsto\lambda_a,
  \]
  where $\lambda_a\colon A\to A$, $b\mapsto -a+a\circ b$,
  is a group homomorphism.
\end{lem}

The proof of the lemma is straightforward and is left as an exercise for the reader.
Direct calculations show that
\begin{align}
  \label{eq:formulas}
   & a\circ b = a+\lambda_a(b),
   &                            & a+b=a\circ \lambda^{-1}_a(b),
   &                            & \lambda_a(a')=-a
\end{align}
hold for $a,b\in A$. Moreover, if
\[
  a*b=\lambda_a(b)-b=-a+a\circ b-b,
\]
then the following identities are easily verified:
\begin{align}
  \label{eq:commutator1} & a*(b+c)=a*b+b+a*c-b,            \\
  \label{eq:commutator2} & (a\circ b)*c=(a*(b*c))+b*c+a*c.
\end{align}

\begin{defn}
  A map $f\colon A\to B$ between two skew braces $A$ and $B$ is said to be a
  \emph{homomorphism} of skew braces if $f(x\circ y)=f(x)\circ f(y)$
  and $f(x+y)=f(x)+f(y)$ for all $x,y\in A$.  The \emph{kernel} of $f$ is
  \[
    \ker f=\{x\in A:f(x)=0\}.
  \]
\end{defn}

\begin{defn}
  Let $A$ be a skew brace. An \emph{ideal} of $A$ is a subset $I$ such that $(I,\circ)$ is a normal subgroup of $(A,\circ)$, $(I,+)$ is a normal subgroup of $(A,+)$ and $\lambda_a(I)\subseteq I$ for all $a\in A$.
\end{defn}

For a group $G$, we write $Z(G)$ to denote the center of $G$.

\begin{exa}
  Let $A$ be a skew brace. Two ideals of $A$ are the \emph{socle} 
  \[
	  \Soc(A)=\ker\lambda\cap Z(A,+)
\]
and the \emph{annihilator} 
\[ 
	\Ann(A)=\Soc(A)\cap Z(A,\circ) 
\]
of $A$. The annihilator was introduced in \cite{MR3917122}.
\end{exa}

The kernel of a skew brace homomorphism is an ideal. If $A$ is a skew brace and $I$ is an ideal of $A$, then $a+I=a\circ I$ holds for all $a\in A$. Then $A/I$ itself forms a skew brace, and the canonical map $A\to A/I$ is a surjective skew brace homomorphism with kernel $I$. This observation aligns with the isomorphism theorems, which hold true for skew braces, similar to their counterparts in ring theory and group theory.

Skew braces exhibit striking similarities with rings, establishing a noteworthy analogy between the two structures. This analogy not only serves as a motivation for various ideas and results but also provides a foundation for exploring and discussing the properties of skew braces. Explorations of this analogy can be found in references such as \cite{MR3177933} and \cite{MR4256133}. Building upon this analogy, a notion of nilpotency in the context of skew braces emerges.

\begin{defn}
  We say that a skew brace $A$ is \emph{right nilpotent} if
  there exists $m$ such that
  $A^{(m)}=\{0\}$, where we define recursively 
  $A^{(1)}=A$ and 
  $A^{(k+1)}=A^{(k)}*A$ for $k\geq1$, where $X*Y$ denotes the additive subgroup of $A$ generated by elements of the form
  $x*y$ for $x\in X$ and $y\in Y$. 
\end{defn}

Later, we will explore the concept of right nilpotency in the context of solutions, and we will discover an insightful interpretation.

There are several problems concerning the algebraic structure of skew braces that we will not discuss here. While several papers have explored the construction of finite simple\footnote{A skew brace is \emph{simple} if it has no non-zero proper ideals.} skew braces of abelian type \cite{MR4020748,MR4161288,MR4122077}, our understanding of these objects is still quite limited. Moreover, almost nothing about simple skew braces is known when the additive group is not  abelian.

\section{Skew braces and radical rings}

Let $R$ be a ring. The
operation
$R\times R\to R$, $(x,y)\mapsto x\circ y=x+xy+y$, is always associative and
$x\circ 0_R=0_R\circ x=x$ for all $x\in X$.

\begin{defn}
  We say that a non-unitary ring $R$ is a \emph{radical ring} if $(R,\circ)$ is a group.
\end{defn}

Several characterizations of radical rings have been developed in the literature \cite{MR3308118}. However, for our purposes, it is enough to work with the notion provided above, which serves as a sufficient working definition.

\begin{exa}
  The subset
  \[
    R=\left\{\frac{2x}{2y+1}:x,y\in\mathbb{Z}\right\}\subseteq\mathbb{Q}
  \]
  is a radical ring with the usual
  addition and
  $(u,v)\mapsto u\circ v=u+uv+v$.
\end{exa}

\begin{exa}
  The additive group $\mathbb{Z}/4$ is a radical
  ring with circle operation  $(x,y)\mapsto x\circ y=x+y+2xy$.
\end{exa}

\begin{exa}
  Rings of strictly upper triangular square matrices are examples of radical rings.
\end{exa}

\begin{defn}
  A skew brace $A$ is said to be \emph{two-sided} if
  \begin{equation}
    \label{eq:right_compatibility}
    (a+b)\circ c=a\circ c-c+b\circ c
  \end{equation}
  holds for all $a,b,c\in A$.
\end{defn}

It is worth noting that in the context of two-sided skew braces, we do not require the 
additive group to be abelian.

If $A$ is a two-sided skew brace, then
\[
  a\circ(-b)=a-a\circ b+a,  \quad
  (-a)\circ b=b-a\circ b+b,
\]
hold for all $a,b\in A$.

\begin{thm}[Rump]
  There exists a bijective correspondence
  between radical rings and two-sided
  skew braces of abelian
  type.
\end{thm}

The theorem was proved in \cite{MR2278047};
the correspondence is given by the
formula $-x+x\circ y-y=xy$.
% If
% $A$ is a two-sided
% skew brace of abelian type,
% then the multiplication $(x,y)\mapsto -x+x\circ y-y$
% turns $A$ into a radical ring. Conversely,
% if $A$ is a radical ring with
% multiplication $(x,y)\mapsto xy$, then
% the circle operation $(x,y)\mapsto x+xy+y$
% turns $A$ into a two-sided
% skew brace of abelian type.

A skew brace is \emph{associative} if the operation 
\[ 
  (x,y)\mapsto
  x*y=\lambda_x(y)-y
\] 
is associative. The following
theorem answers a question of
Cedó, Gateva-Ivanova and Smoktunowicz~\cite{MR3818285}

\begin{thm}[Lau]
  \label{thm:Lau}
  If $(A,+,\circ)$ is an associative skew brace of abelian type,
  then $(A,+,*)$ is a radical ring.
\end{thm}

The proof appears in~\cite{MR4136750}.
Theorem \ref{thm:Lau} was proved independently by
Kinyon (unpublished).

The reference \cite{MR4592323} contains new results regarding radical rings and two-sided skew braces. This work provides valuable insights and advancements in understanding these algebraic structures.

\section{Solutions to the Yang--Baxter equation}

In this section, our attention is directed towards (combinatorial) solutions to the YBE. The initial exploration of this topic can be traced back to the papers \cite{MR1722951,MR2095675,MR2368074,MR2383056,MR1637256,MR1769723,MR1809284,MR2024436}. These papers introduced various algebraic tools and approaches for studying solutions to the YBE. They serve as important references in the field and have paved the way for further developments in the theory of YBE solutions.

\subsection{Combinatorial solutions}

For a set $X$ and an integer
$n\geq1$, we write
$X^n=X\times\cdots\times X$ ($n$-times).

\begin{defn}
  A \emph{solution} to the Yang--Baxter equation (YBE) is a pair $(X,r)$,
  where $X$ is a non-empty set and $r\colon X^2\to X^2$ is a bijective map such that
  \[
    r_1r_2r_1=r_2r_1r_2
  \]
  where, if
  \[
    r(x,y)=(\sigma_x(y),\tau_y(x)),
  \]
  then
  \begin{align*}
     & r_1\colon X^3\to X^3, &  & r_1(x,y,z)=(\sigma_x(y),\tau_y(x),z), \\
     & r_2\colon X^3\to X^3, &  & r_2(x,y,z)=(x,\sigma_y(z),\tau_z(y)).
  \end{align*}
  The solution $(X,r)$ is said to be \emph{finite} if $X$ is a finite set.
\end{defn}

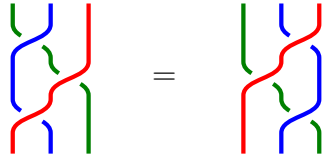
\begin{figure}[H]
  \centering
  \begin{tikzpicture}
    \pic[
      braid/.cd,
      number of strands=3,
      height=.5cm,
      width=.5cm,
      ultra thick,
      gap=0.2,
      name prefix=braid,
      strand 1/.style={red},
      strand 2/.style={blue},
      strand 3/.style={green!50!black},
    ] {braid={a_{1}a_2a_1}};
\node[font=\Large] at (2,1) {\(=\)};
\end{tikzpicture}
\hspace{.5cm}
\begin{tikzpicture}
  \pic[
    braid/.cd,
    number of strands=3,
    height=.5cm,
    width=.5cm,
    ultra thick,
    gap=0.2,
    name prefix=braid,
    strand 1/.style={red},
    strand 2/.style={blue},
    strand 3/.style={green!50!black},
  ] {braid={a_{2}a_1a_2}};
\end{tikzpicture}
\caption{The Yang--Baxter (or braid) equation.}
\label{fig:braid}
\end{figure}

\begin{exa}
  For any non-empty set $X$, the pair $(X,\id_{X^2})$
  is a solution to the YBE.
\end{exa}

Without additional assumptions, the task of finding solutions becomes highly unpredictable. Therefore, we will focus on solutions that meet specific extra assumptions. Given the combinatorial nature of the problem and our intention to utilize group theory, it is compelling to investigate the following intriguing class of solutions: 

\begin{defn}
  We say that a solution $(X,r)$, where 
  \[ 
    r(x,y)=(\sigma_x(y),\tau_y(x)), 
  \] 
  to the YBE 
  is \emph{non-degenerate} if 
  the maps $\sigma_x$ and $\tau_y$ are
  bijective maps for $x,y\in X$.
\end{defn}

With the non-degeneracy assumption, exploring solutions to the YBE becomes even more intriguing as we can now leverage groups that inherently act on our solutions. By incorporating group actions, we gain a richer understanding of the equation's behavior and its connections to various mathematical structures.

\begin{convention}
  A \emph{solution} will always mean a non-degenerate finite solution to the YBE. Moreover,
  if $(X,r)$ is a solution, we will write
  \[
    r(x,y)=(\sigma_x(y),\tau_y(x)).
  \]
\end{convention}

The majority of research on the YBE has primarily focused on non-degenerate solutions. However, there has been a recent surge of interest in studying specific classes of degenerate solutions; see \cite{ACM,CM,CJKVAV,MR4146852}.

\begin{defn}
  A solution $(X,r)$ is said to be \emph{involutive}
  if $r^2=\id_{X^2}$.
\end{defn}

If $(X,r)$ is an involutive solution, then
$r^2(x,y)=(x,y)$ for all $x,y\in X$. Thus
\[
  \tau_y(x)=\sigma^{-1}_{\sigma_x(y)}(x)
\]
for all $x,y\in X$.
This observation implies that in involutive solutions, we only need to know one of the sets $\{\sigma_x: x\in X\}$ or $\{\tau_x: x\in X\}$ to determine the other set. This simplification makes the study of involutive solutions somewhat easier.

\begin{exa}
  \index{Solution!trivial}
  Let $X\ne\emptyset$ be a set. Then $(X,r)$,
  where $r(x,y)=(y,x)$, is an involutive solution. It is called the \emph{trivial solution} over $X$. 
\end{exa}

\begin{exa}
  Let $X\ne\emptyset$ be a set and $\sigma$ and $\tau$ be
  commuting bijections on $X$. Then
  $(X,r)$, where $r(x,y)=(\sigma(y),\tau(x))$, is a solution.
  This is known as the \emph{permutation solution} associated
  with $\sigma$ and $\tau$. This solution is involutive if and only if $\tau^{-1}=\sigma$.
\end{exa}

\begin{exa}
  \label{exa:conjugation}
  Let $X$ be a (union of conjugacy classes of a) group. Then $(X,r)$, where
  $r(x,y)=(y,y^{-1}xy)$
  is a solution.
\end{exa}

\begin{exa}
  \label{exa:core}
  Let $X$ be a group. If 
  $r(x,y)=(xy^{-1}x,x)$, then $(X,r)$ is a solution.
\end{exa}

\begin{exa}
  \label{exa:Alexander}
  Let $X$ be an abelian group and
  $g\in\Aut(X)$. Then $(X,r)$, where
  $r(x,y)=(x-g(x-y), x)$
  is a solution.
\end{exa}

Examples~\ref{exa:conjugation},~\ref{exa:core}, and~\ref{exa:Alexander} exhibit a connection to the theory of racks and quandles; see \cite{MR975077}. These particular solutions play a significant role in constructing knot invariants. See~\cite{MR2896084} for a friendly introduction to the use of racks and quandles in knot theory.

\begin{exa}
  \label{exa:Wada}
  Let $X$ be a group. Then the maps
  \begin{enumerate}
    \item $r(x,y)=(y,x^{-1})$,
    \item $r(x,y)=(y^{-1},x^{-1})$,
    \item $r(x,y)=(x^2y,y^{-1}x^{-1}y)$,
  \end{enumerate}
  satisfy the YBE.
  These solutions were originally discovered by Wada in \cite{MR1167178} and are closely linked to the representation theory of braid groups.
\end{exa}

The following example will be important later.

\begin{exa}
  \label{exa:4_13}
  Let $X=\{1,2,3,4\}$ and $r(x,y)=(\sigma_x(y),\tau_y(x))$ be given~by
  \begin{align*}
     & \sigma_1=(12), &  & \sigma_2=(1324), &  & \sigma_3=(34), &  & \sigma_4=(1423), \\
     & \tau_1=(14),   &  & \tau_2=(1243),   &  & \tau_3=(23),   &  & \tau_4=(1342).
  \end{align*}
  Then $(X,r)$ is an involutive solution.
\end{exa}

\begin{defn}
  Let $(X,r)$ and $(Y,s)$ be solutions. A map $f\colon X\to Y$ is said to be a \emph{homomorphism} of solutions if
  \[
    (f\times f)r=s(f\times f).
  \]
  An \emph{isomorphism} of solutions is
  a bijective homomorphism of solutions.
\end{defn}

Determining the exact number of finite solutions is a complex task. However, in Table~\ref{tab:IYB}, we present the number of non-isomorphic involutive solutions of size up to 10. It is worth noting that for sizes up to 7, the numbers in Table~\ref{tab:IYB} coincide with those reported in~\cite{MR1722951}. However, for size 8, there is a discrepancy of two solutions, as two size-eight solutions are not included in~\cite{MR1722951}.
The numbers in Tables~\ref{tab:IYB} and~\ref{tab:YB} were computed in~\cite{MR4405502}.

\begin{table}[H]
  \centering
  \caption{Enumeration of involutive solutions.}
  \begin{tabular}{|r|ccccccccc|}
    \hline
    $n$ & 2 & 3 & 4  & 5  & 6   & 7    & 8     & 9      & 10\tabularnewline
    \hline
        & 2 & 5 & 23 & 88 & 595 & 3456 & 34530 & 321931 & 4895272\tabularnewline
    \hline
  \end{tabular}
  \label{tab:IYB}
\end{table}

\begin{problem}
\label{pro:small}
How many 
involutive solutions of ``small'' size are there?
\end{problem}

Problem~\ref{pro:small} remains an open question for sizes greater than or equal to 11 in the case of involutive solutions, and for sizes greater than or equal to 9 in the case of non-involutive solutions.

\begin{table}[H]
  \centering
  \caption{Enumeration of non-involutive solutions.}
  \begin{tabular}{|c|ccccccc|}
    \hline
    $n$ & 2 & 3  & 4   & 5    & 6      & 7       & 8\tabularnewline
    \hline
        & 2 & 21 & 253 & 3519 & 100071 & 4602720 & 422449480\tabularnewline
    \hline
  \end{tabular}
  \label{tab:YB}
\end{table}

The following problem is more challenging:

\begin{problem}
Estimate the number of isomorphism classes
of (involutive) solutions of size $n$ for
$n\to\infty$.
\end{problem}

% For $n\geq2$, the \emph{braid group} $\mathbb{B}_n$ is defined as the group with generators $\sigma_1,\dots,\sigma_{n-1}$ and relations
% \begin{align*}
%    & \sigma_i\sigma_{i+1}\sigma_i=\sigma_{i+1}\sigma_i\sigma_{i+1} &  & \text{if }1\leq i\leq n-2, \\
%    & \sigma_i\sigma_j=\sigma_j\sigma_i                             &  & \text{if }|i-j|> 1.
% \end{align*}
% Let $(X,r)$ be a solution.
% For $i<n$ let
% \[
%   r_{i}=\id_{X^{i-1}}\times r\times\id_{X^{n-i-1}}\colon X^n\to X^n.
% \]
% Then the map $\sigma_i\mapsto r_{i}$ extends
% to an action of $\mathbb{B}_n$ on $X^n$.

\subsection{Solutions and skew braces}

A crucial observation is that skew braces yield solutions.

\begin{thm}
  \label{thm:YB}
  Let $A$ be a skew brace. Then
  $(A,r_A)$, where
  \[
    r_A\colon A\times A\to A\times A,\quad
    r_A(x,y)=( -x+x\circ y,(-x+x\circ y)'\circ x\circ y),
  \]
  is a solution to the YBE. Moreover, $(A,r_A)$ is involutive if and only if $A$ is of abelian type.
\end{thm}

Theorem \ref{thm:YB} was proved in \cite[Theorem 3.1]{MR3647970}.
Similar results appear in~\cite{MR1722951,MR1769723,MR2278047,MR1809284}.

The following example highlights the significance and applicability of skew braces in understanding solutions.

\begin{exa}
  The skew brace structure of exact factorizable groups (see Example~\ref{exa:WX}) provides a framework that recovers the solutions discovered by Weinstein and Xu in their work \cite{MR1178147} on factorizable Poisson Lie groups.
\end{exa}

The (minimal) \emph{exponent} $\exp(G)$ of a  finite group $G$ is the least positive integer $n$ such that
$g^n=1$ for all $g\in G$.

\begin{thm}
  \label{thm:|r|}
  Let $A$ be a finite skew brace with more than one
  element and let $G$ be the additive group of $A$.
  If $(A,r_A)$ is the solution to the YBE associated with $A$,
  then the permutstion $r_A$ of $A^2$ has order $2\exp(G/Z(G))$.
\end{thm}

The theorem was proved in~\cite[Theorem 4.13]{MR3763907}.

\begin{defn}
\label{def:structure_group}
  Let $(X,r)$ be a solution.
  The \emph{structure group} of $(X,r)$
  is the group $G(X,r)$ with generators $X$ and relations
  \[
    xy=uv
  \]
  whenever $r(x,y)=(u,v)$.
\end{defn}

Definition~\ref{def:structure_group} goes back to Etingof, Schedler and Soloviev~\cite{MR1722951}. 

In the context of involutive solutions, various properties of their structure groups have been established. For instance, structure groups contain a free abelian group of finite index. Furthermore, it has been determined that these structure groups are Garside groups, as noted in \cite{MR2764830}. This implies that they are torsion-free.
Additionally, these structure groups have been identified as Bieberbach groups, as outlined in \cite{MR1637256}.

\begin{exa}
  If $X\ne\emptyset$ is a set and 
  $(X,r)$ is the trivial solution, then $G(X,r)$ is free
  abelian of rank $|X|$.
\end{exa}

\begin{exa}
  Let $(X,r)$ be the solution of Example~\ref{exa:4_13}.
  The structure group $G(X,r)$ has
  generators
  $x_1,x_2,x_3,x_4$ and relations
  \begin{align*}
     & x_1^2=x_2x_4,
     &                & x_1x_3=x_3x_1,
     &                & x_1x_4=x_4x_3, \\
     & x_2x_1=x_3x_2,
     &                & x_2^2=x_4^2,
     &                & x_3^2=x_4x_2.
  \end{align*}
\end{exa}

\begin{thm}
  \label{thm:iota}
  Let $(X,r)$ be a solution. Then there exists
  a unique skew brace structure over $G(X,r)$ such that $r_{G(X,r)}$ satisfies
  \[
    r_{G(X,r)}(\iota\times\iota)=(\iota\times\iota)r,
  \]
  where $\iota\colon X\to G(X,r)$ is the canonical map. Moreover, if $(X,r)$ is involutive, $G(X,r)$ is a skew brace of abelian type and
  the map $\iota$
  is injective.
\end{thm}

In the non-involutive case, it is important to note that the injectivity of the canonical map is not guaranteed.

\begin{exa}
  Let $X=\{1,2,3,4\}$, $\sigma=(12)$ and $\tau=(34)$. Then $(X,r)$, where $r(x,y)=(\sigma(y),\tau(x))$, is a non-involutive permutation solution. The canonical map
  $\iota\colon X\to G(X,r)$, $i\mapsto x_i$, is not injective, as
  \[
    x_1x_1=x_{\sigma(1)}x_{\tau(1)}=x_2x_1
  \]
  implies that $x_1=x_2$ in $G(X,r)$.
\end{exa}

In general, the additive group of $G(X,r)$ can be shown to be isomorphic to the group $A(X,r)$, which is generated by the elements of $X$ and subject to the relations
\[
  xu=u\sigma_u(v)
\]
whenever $r(x,y)=(u,v)$.

\begin{exa}
  Let $(X,r)$ be an involutive solution. Then $A(X,r)$ is free abelian of rank $|X|$. 
\end{exa}

\begin{exa}
  Let $X=\{1,2,3\}$ and $r(x,y)=(\sigma_x(y),\tau_y(x))$, where
  \[
    \sigma_1=\sigma_2=\sigma_3=(23),\quad
    \tau_1=\id,\quad \tau_2=(132),\quad \tau_3=(123).
  \]
  Then $(X,r)$ is a non-involutive solution
  and $G(X,r)$ has generators $x_1,x_2,x_3$ and
  relations
  \[
    x_1x_2=x_3^2=x_2x_1,\quad
    x_1x_3=x_2^2=x_3x_1,\quad
    x_2x_1=x_1x_2,\quad
    x_3x_1=x_1x_3.
  \]
  Moreover,
  \begin{align*}
    % G(X,r)&=\langle x_1,x_2,x_3:
    % x_1x_2=x_3^2=x_2x_1,
    % x_1x_3=x_2^2=x_3x_1,x_2x_1=x_1x_2,
    % x_3x_1=x_1x_3\rangle,\\  
    A(X,r) & =\langle x_1,x_2,x_3:x_1x_2=x_3x_1=x_2x_3,x_1x_3=x_2x_1=x_3x_2\rangle.
  \end{align*}
\end{exa}

The group $G(X,r)$ naturally acts on $A(X,r)$ via automorphisms, as described in Lemma~\ref{lem:lambda}. Moreover, there exists a bijective 1-cocycle from $G(X,r)$ to $A(X,r)$; see~\cite[Proposition 1.11]{MR3647970}.
 
\begin{exa}
   Let $(X,r)$ be the involutive solution of Example~\ref{exa:4_13}. Then 
   $A(X,r)\simeq\mathbb{Z}^4$. Let $e_1,\dots,e_4$ be the standard basis of $\mathbb{Z}^4$. 
   The map $X\to\mathbb{Z}^{4}$ given by   
   \[ 
    x_i\mapsto e_i,\quad i\in\{1,2,3,4\},
  \]
  %  \[
  %    x_\mapsto \begin{pmatrix}
  %            1\\
  %            0\\
  %            0\\
  %            0
  %    \end{pmatrix},\quad
  %    x_2\mapsto \begin{pmatrix}
  %            0\\
  %            1\\
  %            0\\
  %            0
  %    \end{pmatrix},\quad
  %    x_3\mapsto \begin{pmatrix}
  %            0\\
  %            0\\
  %            1\\
  %            0
  %    \end{pmatrix},\quad
  %    x_4\mapsto \begin{pmatrix}
  %            0\\
  %            0\\
  %            0\\
  %            1
  %    \end{pmatrix},
  %    \]
     extends to a 
     bijective 1-cocycle $G(X,r)\to\mathbb{Z}^4$, where the action of $G(X,r)$ on $\mathbb{Z}^{4}$ is induced by $x_ie_j=e_{\sigma_i(j)}$, $i,j\in X$.  
 \end{exa}

\begin{thm}
  \label{thm:universal}
  Let $(X,r)$ be a solution.
  If $A$ is a skew brace and $f\colon X\to A$ is a map such that
  \[
    (f\times f)r=r_A(f\times f),
  \]
  then there exists a unique homomorphism $\varphi\colon G(X,r)\to A$ of skew braces such that $\varphi\iota = f$ and $(\varphi\times\varphi)r_{G(X,r)}=r_A(\varphi\times\varphi)$.
\end{thm}

Theorems~\ref{thm:iota} and~\ref{thm:universal} appear in \cite{MR3861714,MR3763907}.
Similar results appear in other articles such as~\cite{MR1722951,MR1769723,MR1809284}.

The following result is derived from~\cite{MR1722951}; for more detailed information,
see~\cite[Theorem 2.1]{MR4062375}.

\begin{thm}[Etingof--Schedler--Soloviev]
\label{thm:ESS_representation}
  Let $(X,r)$ be an involutive solution of size $n$.
  Then $G(X,r)$ admits a faithful linear representation inside $\operatorname{GL}_{n+1}(\mathbb{Z})$.
\end{thm}

\begin{exa}
  Let $(X,r)$ be the involutive solution of Example~\ref{exa:4_13}. The structure group admits
  the following faithful linear representation:
  \begin{align*}
     & x_1\mapsto\left(\begin{smallmatrix}
                         0 & 1 & 0 & 0 & 1\\
                         1 & 0 & 0 & 0 & 0\\
                         0 & 0 & 1 & 0 & 0\\
                         0 & 0 & 0 & 1 & 0\\
                         0 & 0 & 0 & 0 & 1
                       \end{smallmatrix}\right),
     &                                     &
    x_2\mapsto\left(\begin{smallmatrix}
                      0 & 0 & 0 & 1 & 0\\
                      0 & 0 & 1 & 0 & 1\\
                      1 & 0 & 0 & 0 & 0\\
                      0 & 1 & 0 & 0 & 0\\
                      0 & 0 & 0 & 0 & 1
                    \end{smallmatrix}\right),
    \\
     & x_3\mapsto\left(\begin{smallmatrix}
                         1 & 0 & 0 & 0 & 0\\
                         0 & 1 & 0 & 0 & 0\\
                         0 & 0 & 0 & 1 & 1\\
                         0 & 0 & 1 & 0 & 0\\
                         0 & 0 & 0 & 0 & 1
                       \end{smallmatrix}\right),
     &                                     &
    x_4\mapsto\left(\begin{smallmatrix}
                      0 & 0 & 1 & 0 & 0\\
                      0 & 0 & 0 & 1 & 0\\
                      0 & 1 & 0 & 0 & 0\\
                      1 & 0 & 0 & 0 & 1\\
                      0 & 0 & 0 & 0 & 1
                    \end{smallmatrix}\right).
  \end{align*}
\end{exa}

It has been established that also structure groups of non-involutive solutions possess a linear structure. However, it is worth noting that the description of these structure groups is not as elegant or straightforward as in the case of involutive solutions.

In the context of Garside theory, Dehornoy made a significant contribution in \cite{MR3374524}. He constructed a faithful linear representation of the structure group of involutive solutions of size $n$ inside $\operatorname{GL}_n(\mathbb{Z}[t,t^{-1}])$.

It is remarkable that algorithmically  all non-degenerate solutions can be derived from skew braces. In the case of involutive solutions, the  algorithm is presented in \cite{MR3527540}. Similarly, for not necessarily involutive solutions, an algorithm is provided in \cite{MR3835326}. These algorithms offer a systematic approach to generating non-degenerate solutions and further strengthen the connection between skew braces and solutions.

\section{Special classes of solutions}

In this section, our focus will be on discussing two distinct families of solutions: indecomposable solutions and multipermutation solutions. We will specifically examine involutive solutions within these families. It is worth noting that our knowledge about these families in the non-involutive case is currently limited, and there is still much to be explored and understood.

\subsection{Indecomposable solutions}

The concept of indecomposable solutions and its corresponding terminology were originally introduced by Etingof, Schedler, and Soloviev in~\cite{MR1722951}.
Indecomposable solutions are intensively studied; see~\cite{MR4565655,MR4502393,MR4568107,MR4308636,MR4514462,MR2132760}.

\begin{defn}
  A solution $(X,r)$ is said to be \emph{decomposable}
  if $X$ can be decomposed into a disjoint union $X=Y\cup Z$ with $Y$ and $Z$ non-empty sets such that $r(Y\times Y)\subseteq Y\times Y$ and 
  $r(Z\times Z)\subseteq Z\times Z$. 
\end{defn}

In this case, $(Y,r|_{Y\times Y})$ 
and $(Z,r|_{Z\times Z})$ are solutions. 

\begin{defn}
  A solution is said to be \emph{indecomposable} if it is not decomposable.
\end{defn}

\begin{exa}
  Let $X\ne\emptyset$ be a set. 
  A permutation involutive solution $(X,r)$, 
  where $r(x,y)=(\sigma(y),\sigma^{-1}(x))$ 
  is indecomposable if and only if $\sigma$ is a cycle of length $|X|$.  
\end{exa}

\begin{exa}
  The solution of Example~\ref{exa:4_13} is indecomposable.
\end{exa}

\begin{exa}
  Let $p$ be a prime number.
  An involutive indecomposable
  solution over a set of cardinality $p$ is isomorphic to $(X,r)$, where
  $X=\mathbb{Z}/p$ and $r(x,y)=(y-1,x+1)$; see~\cite[Theorem 2.13]{MR1722951}. 
\end{exa}

The following family of examples appears in~\cite{MR4514462}.

\begin{exa}
  \label{exa:T}
  Let $(X,r)$ be an involutive solution
  such that the (always invertible) map $T\colon X\to X$,
  $T(x)=\sigma^{-1}_x(x)$, is a cycle of length $|X|$. Then $(X,r)$ is indecomposable.
\end{exa}

Currently, a comprehensive classification of the indecomposable solutions in Example~\ref{exa:T} remains unknown.

\begin{defn}
  Let $(X,r)$ be an involutive solution.
  The group
  \[
    \mathcal{G}(X,r)=\langle\sigma_x:x\in X\rangle
  \]
  is called the \emph{permutation group}
  of $(X,r)$.
\end{defn}

Note that, if $X$ is finite, $\mathcal{G}(X,r)$ is a finite group.
This group acts on the set $X$ by evaluation. Furthermore, it is a skew brace of abelian type. 
% FIXME: add a reference 
Recently, some papers have adopted the theory of skew braces as a means to study the indecomposability of solutions; see~\cite{MR4116644,MR4163866,MR3771874}. 

\begin{pro}
  Let $(X,r)$ be an involutive solution.
  Then $(X,r)$ is indecomposable if and only if
  $\mathcal{G}(X,r)$ acts transitively on $X$.
\end{pro}

Currently, the number of available examples of indecomposable solutions is limited. To gain a deeper understanding of this class of solutions, it is essential to have a greater number of examples at our disposal.

\begin{problem}
\label{pro:indecomposable}
Construct (or enumerate) isomorphism classes of indecomposable solutions of small size.
\end{problem}

The theory of transitive groups of low degree holds potential in tackling Problem~\ref{pro:indecomposable}. It could provide useful insights and approaches towards understanding and solving the problem.

\subsection{Multipermutation solutions}

Multipermutation involutive solutions were introduced 
in~\cite{MR1722951}
for involutive solutions.
For $(X,r)$ an involutive solution,
we define the following relation on the set $X$:
\[
  x\sim y\Longleftrightarrow \sigma_x=\sigma_y,\quad x,y\in X.
\]

A direct calculation shows that $\sim$ is an equivalence relation on $X$.
This equivalence relation induces a solution on the set $X/{\sim}$ of equivalence classes:
\[
  \mathrm{Ret}(X,r)=(X/{\sim},\overline{r}),
  \quad
  \overline{r}([x],[y])=([\sigma_x(y)],[\tau_y(x)]).
\]
the \emph{retraction} of $X$. One defines inductively
\[
  \Ret^1(X,r)=\Ret(X,r),\quad
  \Ret^{k+1}(X,r)=\Ret(\Ret^k(X,r))\text{ for $k\geq1$.}
\]

\begin{defn}
  An involutive solution $(X,r)$ is a
  \emph{multipermutation solution} if there
  exists $n\geq1$ such that $\mathrm{Ret}^n(X,r)$ has only one element.
\end{defn}

Upon examining the list of small involutive solutions, it becomes apparent that a significant majority (over 90\%) of these solutions fall into the category of multipermutation solutions. This observation naturally leads to the question:

\begin{problem}
Can we establish a rigorous proof that ``almost all'' solutions are indeed multipermutation solutions?
\end{problem}

To address this question, it is necessary to provide a precise formulation of the statement, possibly by demonstrating the existence of a certain limit or establishing a quantitative measure of the proportion of multipermutation solutions in the set of all solutions.

\begin{exa}
  \label{exa:5_7}
  Let $X=\{1,2,3,4,5\}$ and
  $r(x,y)=(\sigma_x(y),\tau_y(x))$, where
  \begin{align*}
    &\sigma_1=\sigma_2=\sigma_3=\tau_1=\tau_2=\tau_3=\operatorname{id},&
    &\sigma_4=\tau_4=(45),&
    &\sigma_5=\tau_5= (23)(45).
  \end{align*}
  Then $\Ret(X,r)$ has three elements, $\Ret^2(X,r)$ is the trivial solution 
  over the set of two elements and $\Ret^3(X,r)$ has only one element. Then $(X,r)$ is a multipermutation solution.
\end{exa}

\begin{exa}
  Let $(X,r)$ be the involutive solution of Example \ref{exa:4_13}.
  Then $(X,r)$
  is not a
  multipermutation solution, as
  $\Ret(X,r)=(X,r)$.
\end{exa}

Multipermutation solutions are deeply connected with the theory of skew braces. 
The following result is~\cite[Theorem 5.15]{MR3861714}.

\begin{thm}[Gateva-Ivanova]
  Let $(X,r)$ be an involutive solution. The following statements are equivalent:
  \begin{enumerate}
    \item The solution $(X,r)$ is multipermutation.
    \item The skew brace $G(X,r)$ is right nilpotent.
    \item The skew brace $\mathcal{G}(X,r)$ is right nilpotent.
  \end{enumerate}
\end{thm}

What other ways of detecting multipermutation solutions are there? 
Let $p$ be a prime number and
$G$ be a finite $p$-group. For $k\geq1$,
let
\[
  G^k=\langle g^k:g\in G\rangle.
\]
Then $G^k$ is a normal subgroup of $G$.
We say that $G$ is
\emph{powerful}
if the following conditions hold: if $p>2$,
then $G/G^p$ is abelian; or if $p=2$, then
$G/G^4$ is abelian. The notion goes back to Lubotzky and Mann~\cite{MR0873681}. 

\begin{problem}[Shalev--Smoktunowicz]
  Let $p$ be a prime number. Let $A$ be a skew
  brace of abelian type of size $p^m$ such that 
  $(A,\circ)$ is powerful. Is $A$ 
  right nilpotent? 
\end{problem}

The following problem was posed approximately 80 years ago and remains an unsolved\footnote{In December 2023, Gardam 
proved that the Kaplansky unit conjecture for group rings is false in characteristic zero; see \cite{Gardam}.}
question to this day; see \cite{MR0096696}.

\begin{problem}[Kaplansky's unit problem]
\label{pro:Kaplansky}
Let $G$ be a torsion-free group. Does the group algebra $\mathbb{C}[G]$ have only trivial units?
\end{problem}

Recall that a \emph{trivial unit} of the group algebra $\mathbb{C}[G]$ is an element of the form
$\lambda g$, where $\lambda\in\mathbb{C}\setminus\{0\}$ and $g\in G$.

Gardam proved that the group algebra
of the Promislow
group over the field of two elements
has non-trivial units; see~\cite{MR4334981}.

Problem~\ref{pro:Kaplansky} has an affirmative solution for torsion-free abelian groups. However, the proof of this fact unveils a more general insight: one can substitute the term ``abelian'' with another property and still obtain an affirmative solution.

\begin{defn}
  A group $G$ is said to be \emph{left-orderable}
  if $G$ admits a total order $<$ such that
  $x<y$ implies $zx<zy$ for all $x,y,z\in G$.
\end{defn}

Torsion-free abelian groups, free groups, braid groups
and fundamental group of knots are left-orderable; see~\cite{MR3560661} for other examples.

Left-orderable groups are indeed torsion-free. Moreover, Kaplansky's unit problem has been shown to have an affirmative answer for left-orderable groups. Upon careful examination of the proof, it becomes clear that Problem~\ref{pro:Kaplansky} has an affirmative solution for an even broader class of (torsion-free) groups.

\begin{defn}
  We say that a group $G$ has the \emph{unique product property} if for all
  finite non-empty subsets $A$ and $B$ of $G$ there exists $x\in G$ that can
  be written uniquely as $x = ab$ with $a\in A$ and $b\in B$.
\end{defn}

Left-orderable groups do indeed satisfy the unique product property. Additionally, groups that possess the unique product property are torsion-free. However, it is important to note that the converse does not hold, meaning not all torsion-free groups necessarily have the unique product property.

\begin{exa}[Promislow group]
  Let $P$ be the group
  generated by the matrices
  \begin{align*}
    \begin{pmatrix}
      0 & 1 & 0 & 0   \\
      2 & 0 & 0 & 0   \\
      0 & 0 & 0 & 1/2 \\
      0 & 0 & 1 & 0
    \end{pmatrix}
    \quad\text{ and }\quad
    \begin{pmatrix}
      0 & 0   & 1 & 0 \\
      0 & 0   & 0 & 1 \\
      2 & 0   & 0 & 0 \\
      0 & 1/2 & 0 & 0
    \end{pmatrix}.
  \end{align*}
  Then $P$ is torsion-free and admits the presentation
  \[
    \langle x,y : x^{-1}y^2x=y^{-2}, y^{-1}x^2y=x^{-2} \rangle;
  \]
  see~\cite{MR0470211}. Let
  \begin{multline}
    \label{eq:Promislow}
    S=\{ x^2y,
    y^2x,
    xyx^{-1},
    (y^2x)^{-1},
    (xy)^{-2},
    y,
    (xy)^2x,
    (xy)^2,\\
    (xyx)^{-1},
    yxy,
    y^{-1},
    x,
    xyx,
    x^{-1}
    \}.
  \end{multline}
  A direct tedious calculation shows that each $s\in S^2=\{s_1s_2:s_1,s_2\in S\}$ admits at least two different decompositions
  of the form $s=ab=uv$ for $a,b,u,v\in S$. Thus $P$ does not have the unique product property; see~\cite{MR940281}.
\end{exa}

The following result answers a question of Gateva--Ivanova~\cite{MR3861714}.

\begin{thm}
  \label{thm:BCV}
  Let $(X,r)$ be an involutive solution. Then
  $(X,r)$ is a multipermutation solution if and only if the structure group $G(X,r)$ is left-orderable.
\end{thm}

The implication $\implies$ was proven by Jespers and Okniński \cite{MR2189580} and Chouraqui \cite{MR3572046}.
On the other hand, the reverse implication is notably harder and relies on the theory of skew braces. The details and proof of this implication can be found in \cite{MR3815290}.

\begin{cor}
  Let $(X,r)$ be an involutive solution. If all Sylow subgroups of $\mathcal{G}(X,r)$ are cyclic, then $(X,r)$ is a multipermutation solution. 
\end{cor}

The aforementioned corollary, which can be found in \cite{MR4062375}, makes use of diffuse groups, specifically relying on~\cite[Theorem 3.15]{MR3548136} and~\cite[Theorem 7.12]{MR3974961}. Furthermore, it extends a result initially established in~\cite[Proposition 3.7]{MR1722951} for cyclic groups.

\begin{problem}
\label{pro:UPP}
Let $(X,r)$ be an involutive solution. Under what conditions does the associated group $G(X,r)$ possess the unique product property?
\end{problem}

Problem~\ref{pro:UPP} was presented in~\cite{MR3974961}. To shed light on its importance, there exist several insightful observations worth considering.

There are only few involutive solutions of size up to ten that are not multipermutation solutions. In most cases, it is possible to utilize the set~\eqref{eq:Promislow} to demonstrate that the structure group of these solutions does not possess the unique product property. In fact, in these cases, it can be proven that the structure groups contain a subgroup isomorphic to the Promislow group; see~\cite{MR4062375}.

\begin{exa}
  Let $(X,r)$ be the solution of Example~\ref{exa:4_13}. Then $G(X,r)$ does not have the unique product property. For
  $x=x_1x_2^{-1}$ and $y=x_1x_3^{-1}$
  consider the Promislow set $S$ of~\eqref{eq:Promislow}.
  Then each $s\in S^2$ admits at least two different decompositions
  of the form $s=ab=uv$ for $a,b,u,v\in S$.
\end{exa}

There exist a few non-multipermutation solutions where the structure group does not include the Promislow subgroup as a subgroup. In particular, the Promislow set~\eqref{eq:Promislow} cannot be utilized to prove that the group lacks the unique product property.

\begin{exa}
  \label{exa:candidate}
  Let $X=\{1,\dots,8\}$ and $r(x,y)=(\sigma_x(y),\tau_y(x))$, where
  \begin{align*}
     & \sigma_1=\sigma_2=(3745),     &  & \tau_1=\tau_2=(3648),     \\
     & \sigma_3=\sigma_4=(1826),     &  & \tau_3=\tau_4=(1527),     \\
     & \sigma_5=\sigma_7=(13872465), &  & \tau_5=\tau_7=(16542873), \\
     & \sigma_6=\sigma_8=(17842563), &  & \tau_6=\tau_8=(13562478).
  \end{align*}
  Then $(X,r)$ is not a multipermutation solution, as
  $\Ret(X,r)$ is isomorphic to the solution
  of Example~\ref{exa:4_13}.
  Thus $G(X,r)$ is not left-orderable. Moreover,
  $G(X,r)$
  does not contain subgroups isomorphic to the Promislow subgroup. Note that $G(X,r)$ 
  admits a faithful linear representation given by Theorem \ref{thm:ESS_representation}. 
  Does $G(X,r)$ have the unique product property?
\end{exa}

The structure group of the solution provided in Example~\ref{exa:candidate} serves as a promising candidate for studying Kaplansky's unit problem.

\section{Group theory}

\subsection{Solvable groups}

There are some problems within the theory of skew braces that are directly linked to the solvability of groups. To begin, let us recall an important result that was proven in \cite[Theorem 2.15]{MR1722951}.

\begin{thm}[Etingof--Schedler--Soloviev]
  \label{thm:ESS_solvable}
  Let $A$ be a finite skew brace of abelian type. Then $(A,\circ)$ is solvable.
\end{thm}

It is worth noting that in \cite{MR1722951}, the result was formulated using a different terminology, as the concept of skew braces of abelian type was introduced later.

Theorem \ref{thm:ESS_solvable} raises an intriguing question: Is every finite solvable group the multiplicative group of a skew brace of abelian type? This question has been explored in several papers, including~\cite{MR2584610,MR4198084}. In particular, in \cite{MR3465351}, Bachiller presented a solvable group that does not arise as the multiplicative group of a skew brace of abelian type. This demonstrates that the answer to the aforementioned question is negative.

\begin{problem}
\label{pro:minimal}
Which is the minimal size of a finite solvable group that does not arise as the multiplicative group of a skew brace of abelian type?
\end{problem}

Problem~\ref{pro:minimal} can be regarded as a discrete analog of a conjecture initially posed by Milnor in the context of flat manifolds~\cite{MR0454886}. It is crucial to emphasize that Milnor's conjecture has been disproven, as demonstrated in~\cite{MR1316552,MR1411303}.

\begin{problem}
Which finite solvable groups arise as the multiplicative group of a skew brace of abelian type?
\end{problem}

The problem we will consider next originates from the theory of Hopf-Galois structures and was initially formulated by Byott in \cite{MR3425626}. It is also mentioned as Problem 19.91 in \cite{MR3981599}.

\begin{problem}[Byott]
Let $A$ be a finite skew left brace of solvable type. Is the
multiplicative group of $A$ solvable?
\end{problem}

The answer to Byott's problem is currently unknown and continues to be the subject of active research and investigation. Some partial results related to this problem can be found in \cite{B,MR4210997}.

\subsection{Growth series}

The following problem holds particular interest in the non-involutive case.

The structure skew brace $G(X,r)$ associated with a finite solution $(X,r)$ provides valuable information about the solution itself. Note that the underlying multiplicative group of the skew brace $G(X,r)$ is a finitely generated group.
To further visualize this group structure, we can consider the Cayley graph $\Gamma(X,r)$ of the multiplicative group of $G(X,r)$, using the generating set $X$. In this graph, the vertices correspond to the elements of the multiplicative group of $G(X,r)$, and for every element $g$ in the group and every generator $x$ in $X$, there exists an edge labeled with $x$ connecting $g$ and $gx$ in the graph.

We can define a distance between vertices in the graph.
For $n\geq0$, let
\begin{align*}
   & B(1,n)=\{g\in G(X,r):\operatorname{dist}(1,g)\leq n\},
   &                                                        & \gamma_{G(X,r)}(n)=|B(1,n)|.
\end{align*}
It is important to note that
$\gamma_{G(X,r)}(n)<\infty$.
We say that the solution
$(X,r)$ has \emph{rational growth} if
the \emph{growth series}
\[
  \sum_{n\geq0}\gamma_{G(X,r)}(n)t^n\in\mathbb{Z}[\![t]\!]
\]
of the solution $(X,r)$
is a rational function, that is the quotient
$p(t)/q(t)$ of polynomials $p(t)$ and $q(t)$.

\begin{problem}
Compute the growth series of finite solutions.
\end{problem}

In the work by Benson~\cite{MR714092}, it was established that groups containing a finite-index subgroup isomorphic to $\mathbb{Z}^n$ possess rational growth series. The structure group of a finite solution exhibits a finite-index subgroup that is isomorphic to $\mathbb{Z}^n$. As a consequence, we can conclude that the structure group of a finite solution has a rational growth series.

\section*{Acknowledgements}

This work was partially supported by 
the project OZR3762 of Vrije Universiteit Brussel and 
FWO Senior Research Project G004124N.

\bibliographystyle{abbrv}
\bibliography{refs}

\end{document}